\documentclass[12pt]{article}
\usepackage{amsmath,epsfig,amssymb,float,graphicx,amsfonts, amsbsy}
\usepackage{xspace,colortbl}

\title{Covariance of centered distributions on manifold}
 
\author{\Large Nikolay H. Balov
\footnote{Florida State University, Department of Statistics}\\
{\small\texttt{balov@stat.fsu.edu}}}

\begin{document}

\maketitle

\begin{abstract}

We define and study a family of distributions with domain 
complete Riemannian manifold. They are obtained by projection 
onto a fixed tangent space via the inverse exponential map. 
This construction is a popular choice in the literature for it makes it easy to generalize 
well known multivariate Euclidean distributions. 
However, most of the available solutions use coordinate specific definition that makes them less versatile.
We define the distributions of interest in coordinate independent way by utilizing co-variant 2-tensors.
Then we study the relation of these distributions to their Euclidean counterparts. 
In particular, we are interested in relating the covariance to the tensor 
that controls distribution concentration.
We find approximating expression for this relation in general 
and give more precise formulas in case of manifolds of constant curvature, 
positive or negative.
Results are confirmed by simulation studies of the standard normal distribution 
on the unit-sphere and hyperbolic plane.

\end{abstract}


{\section{Introduction}}

We are interested in defining and studying some of the properties of distributions on 
complete Riemannian manifolds. A typical example of such manifolds is the unit n-sphere $\mathbb{S}^n$.
In this sense, the subject of our study has as a primary application, but not limited to, 
directional statistics, a branch of statistics dealing with directions and rotations in $\mathbb{R}^n$.

Pioneers in the field are Fisher, R.A., \cite{fisher_sphere} and von Mises.
In recent years directional statistics proved to be useful in variety of disciplines like 
shape analysis \cite{mardia-dirs}, geology, crystallography \cite{krieger-crystal}, 
bioinformatics \cite{mardia-protein} and data mining \cite{bahlmann-hand}.

The best known distribution from the field of directional statistics is 
the von Mises-Fisher distribution. It is defined on the unit n-sphere by the density 
$$
f_n({\bf x};\mu, k) = C_n(k) \exp(k\mu' {\bf x}),\\ {\bf x}\in \mathbb{S}^n,
$$
where $k\ge 0$, $\mu\in \mathbb{S}^n$ and normalizing constant $C_n(k)$. It is applied initially 
for studying electric fields (n=2). Its one dimensional variant, the von Mises distribution, 
is also known as the circular normal distribution. 

Another important distribution is Fisher-Bingham-Kent(FBK) distribution, proposed by Kent, J. in 1982.
It is defined on $\mathbb{S}^2$ by the density 
$$
f(\mathbf{x})=\frac{1}{\textrm{c}(\kappa,\beta)}\exp\{\kappa\boldsymbol{\gamma}_{1}\cdot\mathbf{x}+\beta[(\boldsymbol{\gamma}_{2}\cdot\mathbf{x})^{2}-(\boldsymbol{\gamma}_{3}\cdot\mathbf{x})^{2}]\},
$$
where $\gamma_1$, $\gamma_2$ and $\gamma_3$ are three orhonormal disrections in $\mathbb{R}^3$.
A recent application of Kent distribution can be found in \cite{kent-proitein}.

The family of centered distributions we are going to consider 
includes von Mises-Fisher distributions but not FBK distributions which are of mixed nature.
Centered distributions are obtained by projecting the distribution domain onto 
a fixed vector space, namely a tangent space on manifold. This approach is well known and 
easy to implement. However, we think that not all of its aspects are treated rigorously. 
One problem that needs care is defining distributions in coordinate free manner. 
This issue is important when the domain is a compact Riemannian manifold as $\mathbb{S}^n$ and does not 
accept a global parametrization. Another problem arises in the study of covariance, 
which has coordinate specific nature. Only those properties of distributions that are 
coordinate system invariant are relevant in comparison studies.

Here we do not target a specific application, 
but rather aim at generalization and pedagogical improvement over the existing solutions like 
providing coordinate free definition of large class of distributions on complete manifolds. 

Another direction in this study is the impact of domain curvature on the covariance of distributions of interest. 
Again, we improve upon some existing results \cite{pennec-riemann}, by generalizing and being more precise.
Finally, we provide simulation results, something that up to our knowledge is missing in the literature, 
that illustrate and confirm the formal developments on specific spaces of constant curvature, 
the unit 2-sphere and hyperbolic plane.

{\section{Definition of centered distributions}}

Let M be a Riemannian n-manifold, $q\in M$ and 
let $Exp_q$ be the exponential map at $q$, $Exp_q: T_qM \to M$.
If M is complete, then the exponential map $Exp_q$ is defined on the whole tangent space $T_qM$.
Throughout this paper we will assume that M is a complete Riemannian n-manifold. 

There is a maximal open set $B(q)$ in $T_pM$ containing the origin, where $Exp_q$ is a diffeomorphism. 
Then the set $\mathcal{B}(q) = Exp_q(B(q))$ is called maximal normal neighborhood of $q$. 
On this normal neighborhood the exponential map is invertible and let 
$$
Log_q = Exp_q^{-1} : \mathcal{B}(q) \to T_pM
$$
be its inverse, the so called log-map. $Log_q$ is diffeomorphism on $\mathcal{B}(q)$.

The Borel sets on M generated by the open sets on M form a $\sigma$-algebra $\mathcal{A}$ on M. 
Any Riemannian manifold has a natural measure $\mathcal{V}$ on $\mathcal{A}$, 
called {\it volume measure}. In local coordinates $x$ it is given by
$$
dV(x) = \sqrt{|G(x)|}dx,
$$
where $G(x)$ is the matrix representation of the metric tensor, $|G|$ is its determinant and 
$dx$ is the Lebesgue measure in $\mathbb{R}^n$. More details one can find in \cite{chavel-riemannian}, ch. 3.3.

We consider a family $\mathcal{Q}$ of 
distributions on M given by density differentials 
\begin{equation}\label{eq:distrib_tensor_density}
dQ(p; q,T,f) = k f(T(Log_qp, Log_qp)) dV(p),
\end{equation}
where $q\in M$, $T$ is a symmetric and positive definite co-variant 2-tensor (bi-linear form)
at tangent space $T_qM$, 
$f:\mathbb{R} \to \mathbb{R}^+$ is a function on M and k is a normalizing constant.
We call the elements of $\mathcal{Q}$ {\it centered distributions} for an obvious reason - 
their densities are defined via projection onto a single tangent space ($T_qM$) placed at a central point ($q$). 
Note that their intrinsic means may or may not coincide with $q$.
Also, as defined the distributions from $\mathcal{Q}$ are absolute continuous with respect to the volume measure 
with $kf(T(.,.))$ being their densities.

A particular member of the family $\mathcal{Q}$ takes 
$$
T(X,Y) = <X,Y>, \\ X, Y \in T_q M,
$$
$$
f(t) = \exp(-\frac{1}{2\sigma^2}t), 
$$
and defines the so called {\it standard normal distribution} on M at $q$.

Sometimes we want the log-map to be surjective on the entire tangent space at $q$ 
except eventually a subset of measure zero. 
With the current definition of the log-map we have $Log_q(\mathcal{B}(q)) = B(q)$ and when the 
cut locus $Cut(q)$ of $q$ is non-empty, $B(q)$ is a bounded star-like neighborhood in $T_pM$.
Can we extend the definition of $Log_q$ so that it covers the maximal possible image of $\mathcal{B}(q)$, 
$T_qM$? We are going to introduce a multi-value version of the log-map designed 
to meet this requirement.

The set of critical points of $Exp_q$, i.e. the set where $Exp_q$ is not diffeomorphism, 
is closed and with volume measure zero. 
(for more details see \cite{chavel-riemannian}, Th 3.2 and Prop. 3.1). 
In fact, the set of non-critical points of $Exp_q$ is exactly $\mathcal{B}(q)$, the 
maximal normal neighborhood of $q$.
Thus, we have that $\mathcal{B}(q)$ is open in M and $\mathcal{V}(M\backslash\mathcal{B}(q)) = 0$.
For any $p\in\mathcal{B}(q)$, there exists a neighborhood V of p 
such that $V\subset \mathcal{B}(q)$. 
Since $W = Exp_q^{-1}(V)$ is open in $T_pM$, which has a countable basis, 
W has countably many connected components, $W = \cup_{i\ge 1} W_i$.
Moreover, each connected component $W_i$ of $W$ 
maps diffeomorphically on V by $Exp_q$. 
Therefore, if we consider 
$\mathcal{B}(q)$ to be a submanifold of M, then the map
$$
Exp_q: Exp_q^{-1}(\mathcal{B}(q)) \to \mathcal{B}(q),
$$
is a covering of $\mathcal{B}(q)$.
In fact, we can take $V = \mathcal{B}(q)$ and then 
$$
Exp_q |_{W_i}: W_i \to \mathcal{B}(q)
$$
are diffeomorphisms.

Define $Log_q|_{W_i}(p) = v_i$, for the unique $v_i\in W_i$ such that $Exp_q(v_i)=p$. 
The diffeomorphisms $Log_q |_{W_i}:\mathcal{B}(q) \to W_i$ we call {\it leafs} of the log-map.
The multi-value version $\widetilde{Log_q}$ of $Log_q$ is defined on the entire $\mathcal{B}(q)$ by
\begin{equation}\label{eq:logmap_multival}
\widetilde{Log_qp} = \{Log_q|_{W_i}(p)\}_{i\ge 1},\\ p\in \mathcal{B}(q).
\end{equation}
%

We define 
\begin{equation}\label{eq:wrapped_pdf}
f(T(\widetilde{Log_{q}p}, \widetilde{Log_qp})) = \sum_{i=1}^{\infty} f(T(Log_qp|_{W_i}, Log_qp|_{W_i})).
\end{equation}
and then the distribution form (\ref{eq:distrib_tensor_density}) has to be read as
$$
dQ(p; q,T,f) = k f(T(\widetilde{Log_{q}p}, \widetilde{Log_qp})) dV(p) = 
$$
$$
k \sum_{i=1}^{\infty} f(T(Log_qp|_{W_i}, Log_qp|_{W_i})) dV(p).
$$
We refer to the operation (\ref{eq:wrapped_pdf}) as {\it folding} a density. 
Basically, the support of $f$ determines how many leafs of the log-map we use.

Recall that if $(x_1,...,x_n)$ is an orthonormal basis of $T_qM$, 
the normal coordinate system on $W_i$ is given by
$$
v=(v^1,...,v^n) \mapsto \phi(v) = Exp_q(v^1x_1+...+v^nx_n), \\ v^1x_1+...+v^nx_n \in W_i
$$ and  
the log-map is particularly simple, $Log_q|_{W_i}(v) \equiv Log_q \phi(v) |_{W_i} = v$. 

The benefit of introducing (\ref{eq:logmap_multival}) and (\ref{eq:wrapped_pdf}) 
is clear when one integrates functions of the log-map. For example, 
expectation with respect to $Q\in\mathcal{Q}$ of any measurable function $h(p)$ on M is 
$$
E h = 
k \sum_i \int_M h(p) f(T(Log_qp|_{W_i}, Log_qp|_{W_i})) dV(p) = 
$$
$$
k \sum_i \int_{\phi^{-1}( W_i)} h(\phi(v)) f(T(v, v)) dV(\phi(v)) = 
$$
$$
\int_{\mathbb{R}^n} h(\phi(v)) f(T(v, v)) dV(\phi(v)),
$$
where the last integral is the Lebesgue one on the whole $\mathbb{R}^n$.
Using the multi-value log-map all density functions $f$ in $\mathbb{R}^n$, like the normal ones, 
can be manipulated easier on a general manifold M, because we do not change their support.

\newtheorem{densapprox_example_1}{Example}
\begin{densapprox_example_1}
Let $M$, be the unit n-sphere $\mathbb{S}^n$.
Fix a point $q\in\mathbb{S}^n$. 
The cut locus point for q is $-q$, the antipodal point.
Thus, $\mathcal{B}(q) = \mathbb{S}^n\backslash \{-q\}$.
Define $U_k = B_{k\pi}(q)$, the ball on $T_p \mathbb{S}^n$ with radius $k\pi$, $k\ge 1$.
The maximal normal neighborhood for q is $U_1$.
We have $Exp_q^{-1}(\mathcal{B}(q)) = \cup_i W_i$ for $W_i = B_{(i+1)\pi}(q)\backslash\overline{B_{i\pi}(q)}$.
Let $n_{qp}=Log_qp/||Log_qp||$ be the unit tangent vector at q in the direction of p, then 
$Log_q|_{B_1}(p) = d(q,p)n_{qp}$ for 
$$
d(q,p)=cos^{-1}<q,p> \\ \in [0,\pi] 
$$
and 
$$
\widetilde{Log_qp} = \{(d(q,p) \pm 2\pi i)n_{qp} \}_{i\ge 0}.
$$
\end{densapprox_example_1}
\newtheorem{densapprox_remark_1}{Remark}
\begin{densapprox_remark_1}
In a sense, the proposed extension of the log-map with corresponding modified distributions (\ref{eq:wrapped_pdf}) 
is generalization of the concept of wrapped distributions. 
These are densities $f$ on the line, 'wrapped' around the circumference of the unit circle $\mathbb{S}^1$: 
$
f(\theta) = \sum_{i=-\infty}^{\infty} f(\theta + 2\pi i), \theta\in [0,2\pi),
$
as used in \cite{bahlmann-hand}.
\end{densapprox_remark_1}
\newtheorem{densapprox_example_2}[densapprox_example_1]{Example}
\begin{densapprox_example_2}
Von Mises-Fisher distribution is a centered distribution of form (\ref{eq:distrib_tensor_density}) if we take 
$q = \mu$, $T(v,v)=v'v = ||v||^2$, $f(t) = c_0 \exp(k cos(t))$, $t\in[0,2\pi]$ and a normalizing constant $c_0$.
Support of $f$ is bounded and we use only the first leaf of the log-map.
\end{densapprox_example_2}
\newtheorem{densapprox_example_3}[densapprox_example_1]{Example}
\begin{densapprox_example_3}
Gamma distribution on M can be defined by
$$
f(t) = c_0 t^{k-1}\exp(-t/\theta),\\ t\ge 0
$$
for $\theta > 0$ and $T(v,v) = v'v$. Constant $c_0$ is determined by
$$
c_0^{-1} = \int_{\mathbb{R}^n} |v|^{k-1}\exp(-|v|/\theta) dv = 
$$
$$
\int_0^{\infty}(\int_{\mathbb{S}_r^n} d\phi) r^{k-1}\exp(-r) dr  = 
2\pi^{(n+1)/2}\theta^{n+k}\frac{\Gamma(n+k)}{\Gamma((n+1)/2)},
$$
where we used that the area of $\mathbb{S}_r^n$ is $\frac{2\pi^{(n+1)/2}r^n}{\Gamma((n+1)/2)}$.
Because the support of $f$ is the whole $\mathbb{R}$, we have a folded density.
\end{densapprox_example_3}

Unfortunately, both von Mises-Fisher and Gamma multivariate distributions do not have explicit 
expression for their second moments which make them less useful in the context of the following results.

{\section{Approximating the covariance}}

Let $Q$ be a distribution from $\mathcal{Q}$.
Covariance of $Q$ we call 
a contra-variant 2-tensor at tangent space $T_q M$ given by
$$ 
\Sigma = k \int_p (Log_qp) (Log_qp)' f(T(Log_qp, Log_qp)) dV(p).
$$
Note that when $q$ is the mean (intrinsic) of $Q$, $\Sigma$ is a covariance in the usual sense, 
but here we do not require $q$ to be a mean and we use the term covariance in a different context, namely, 
as a quantity measuring the dispersion about the center $q$ of $Q$.

We want to obtain an approximating expression for the covariance of $\mathcal{Q}$ 
as a function of the tensor T and the first few moments of f.

In normal coordinates $v$ at $q$, the volume measure can be approximated by
\begin{equation}\label{eq:approx_volume_manifold}
dV(v) = [1 - \frac{1}{6}v'(Ric)v + O(|v|^3)]dv,
\end{equation}
where $Ric$ is the matrix representation of the Ricci tensor (see for example Th. 2.17 in Chavel \cite{chavel-riemannian}).

X. Pennec \cite{pennec-riemann} used equation (\ref{eq:approx_volume_manifold}) to approximate the covariance of 
normal distribution. His approximation is $\Sigma \thickapprox T^{-1} - \frac{1}{3}T^{-1}(Ric)T^{-1}$.

We use the above approximation of the volume form to obtain more general result applied for 
densities of centered distributions given by (\ref{eq:distrib_tensor_density}). In addition,
we derive more precise variance estimation on the unit 2-sphere and the hyperplane. 
Finally, we provide some simulation results to confirm the formulas.

Let $T$ be the matrix representation of tensor T with respect to coordinates $v$. 
Let $T^{-1} = U\Lambda U'$ be the eigenvalue decomposition of $T^{-1}$ 
with diagonal matrix of eigenvalues $\Lambda$. 
Define $S=U\Lambda^{1/2}$. Then $T^{-1} = SS'$. 
The determinant of $S$ is $|S| = |T|^{-1/2}$ and its norm is $||S||$ given as 
$||S|| = \sup \{||Sx||_2,\\ ||x||_2=1 \}$. 
$||S||$ is the maximal eigenvalue of $S$, which is strictly positive.
Moreover $||S|| \le ||U||||\Lambda^{1/2}|| \le ||T^{-1}||^{1/2} = \lambda_{min}^{-1/2}$, 
where $\lambda_{min}$ is the minimal eigenvalue of $T$.

We change the variables $v$ to $w=(w_i)$ according to 
$$
v = Sw.
$$ 
Then $v'Tv = w'w$ and $vv' = S(ww')S'$.
Density $f$ is assumed to satisfy
\begin{equation}\label{eq:density_assumptions1}
\int_{\mathbb{R}^n} f(w'w)dw = 1,\\ 
\int_{\mathbb{R}^n} wf(w'w)dw = 0,\newline
\end{equation}
and let 
\begin{equation}\label{eq:density_assumptions2}
\int_{\mathbb{R}^n} ww'f(w'w)dw = C,\\  
\int_{\mathbb{R}^n} [(ww')\otimes (ww')] f(w'w)dw = D.
\end{equation}
$C$ is a symmetric and positive definite $n\times n$ matrix, while 
$D$ is the expectation of the Kronecker product $(ww')\otimes (ww')$ and thus, it is a $n^2\times n^2$ matrix.
Let $D = \{D_{kl}^{ij}\}_{klij}$ and for every $k,l\in \{1,...,n\}$, $D_{kl}$ is the corresponding 
$n\times n$ matrix.
Let $R=S'(Ric)S = (r_{ij})$. 
By $tr(RD)$ we will understand the $n\times n$ matrix with elements 
$[tr(RD)]_{ij} = \sum_{k,l} r_{kl} D_{kl}^{ij}$.

Now we are ready to formulate the following
\newtheorem{lemma_approx_density}{Lemma}
\begin{lemma_approx_density}
Under the assumptions (\ref{eq:density_assumptions1}) and (\ref{eq:density_assumptions2}), the density form 
(\ref{eq:distrib_tensor_density}) has normalizing constant
\begin{equation}\label{eq:approx_density_k}
k^{-1} = |S|(1-\frac{1}{6}tr(RC) + \epsilon)
\end{equation}
and covariance 
\begin{equation}\label{eq:approx_density_cov}
k^{-1}\Sigma = |S| S (C - \frac{1}{6}tr(RD) + \epsilon I_n) S'.
\end{equation}
where the function $\epsilon(S) = O(||S||^{3})$.
\end{lemma_approx_density}
{\it The Proof} is a straightforward derivation. First observe that by definition 
\begin{equation}\label{eq:int_def_kinv}
k^{-1} = \int_{\mathbb{R}^n} f(v'Tv)dV(v)
\end{equation}
and 
\begin{equation}\label{eq:int_def_kinvsigma}
k^{-1}\Sigma = \int_{\mathbb{R}^n} vv'f(v'Tv)dV(v).
\end{equation}
assuming $Log_q \mathcal{B}(q) \cong \mathbb{R}^n$, which we can always guarantee 
by folding, eventually, the original density $f$ (see definition (\ref{eq:wrapped_pdf})).
In the rest of this section all integrals are assumed with domain $\mathbb{R}^n$.

We proceed by expressing the terms that appear above when the volume form is replaced by approximation (\ref{eq:approx_volume_manifold}).
Obviously, $\int f(v'Tv)dv = |S|$ and then 
\begin{equation}\label{eq:int_vpv}
\int v'(Ric)v f(v'Tv)dv = \int tr(w'Rw) f(w'w)|S|dw = 
|S| tr(RC).
\end{equation}
Similarly
\begin{equation}\label{eq:int_vvp}
\int vv'f({v'Tv})dv = S (\int ww' f(w'w)|S|dw) S' =  
|S|SC S'.
\end{equation}
Then we derive 
$$
\int vv'(v'(Ric)v) f(v'Tv)dv = 
|S| S (\int ww' (w'Rw)f(w'w)dw) S'
$$
with the $(ij)^{th}$ element of the last integral equal 
$$
\int w_iw_j (w'Rw)f(w'w)dw = 
\int w_iw_j \{ \sum_{k,l} w_k r_{kl} w_l \} f(w'w) dw = 
$$
$$
\sum_{k,l} r_{kl} \int w_iw_j w_k w_l f(w'w)dw = [tr(RD)]_{ij}.
$$
Thus, 
\begin{equation}\label{eq:int_vvp_vpv}
\int vv'(v'(Ric)v) f(v'Tv)dv = |S| Str(RD)S'.
\end{equation}
Finally for the error term we have  
$$
\int |v|^3 f(v'Tv)dv \le \int ||s||^3 |w|^3 |S|dw \le ||S||^{3+n} \int |w|^3 dw,
$$
using the fact that $|S| \le ||S||^{n}$. Since given the assumptions we made the last integral is bounded, we have
\begin{equation}\label{eq:int_v3}
\int |v|^3 f(v'Tv)dv = |S| O(||S||^{3}).
\end{equation}
Plugging (\ref{eq:int_vpv}), (\ref{eq:int_vvp}), (\ref{eq:int_vvp_vpv}) and (\ref{eq:int_v3}) into 
(\ref{eq:int_def_kinv}) and (\ref{eq:int_def_kinvsigma}) 
one obtains the claim.
$\Box$

Formulas (\ref{eq:int_def_kinv}) and (\ref{eq:int_def_kinvsigma}) are given 
with respect to a normal coordinates $v$, which are not unique.
We will show how they change with a change of coordinates and what 
is invariant to such a change.

Let $\tilde{v}$ be another normal coordinate system at $q$ and 
matrix $A$ be the Jacobian of the change from $v$ to $\tilde{v}$, i.e. $\tilde{v}=Av$. 
$A$ is orthogonal matrix, $A\in O(n)$.

Since T is a symmetric positive definite co-variant 2-tensor then
$T^{-1}$ is a contra-variant 2-tensor and so it is $\Sigma$. 
Under the coordinate change we have
$$
T^{-1} \mapsto A T^{-1} A', \\ S \mapsto A S, \textrm{ and } \Sigma \mapsto A \Sigma A'.
$$
Matrices $C$ and $D$ remains unchanged and so does $R=S'(Ric)S$, 
because $Ric$ is a co-variant tensor such that $Ric \mapsto (A^{-1})'(Ric)A^{-1}$.
Moreover 
$$
S^{-1} \Sigma (S^{-1})' \mapsto S^{-1}A^{-1}  A\Sigma A' (A^{-1})'(S^{-1})'
$$
and hence, the above quantity is also coordinate system invariant.
We showed the following

\newtheorem{lemma_approx_density2}[lemma_approx_density]{Lemma}
\begin{lemma_approx_density2}
Matrix $S^{-1} \Sigma (S^{-1})'$ is an invariant to the normal coordinate system at $q$ 
and satisfies
\begin{equation}\label{eq:int_vvp_vpv}
S^{-1} \Sigma (S^{-1})' = \frac{C - \frac{1}{6}tr(RD) + \epsilon I_n}{1 - \frac{1}{6}tr(RC) + \epsilon},
\end{equation}
where $\epsilon(T) = O(||T^{-1}||^{3/2})$.
\end{lemma_approx_density2}

\newtheorem{densapprox_example_4}[densapprox_example_1]{Example}
\begin{densapprox_example_4}
We take a normal distribution on M, defined by 
$$f(v) = (2\pi)^{-n/2}\exp(-\frac{1}{2}v'Tv)$$ for 
a co-variant tensor $T$. Since 
$\int w_i^2f(w'w)dw = 1$, $\int w_i^4f(w'w)dw = 3$ and 
$[tr(RD)]_{ij} = r_{ij} + r_{ji} + n\delta_{ij}r_{ij}$, we have  
$$
C = I_n,\\ tr(RC) = tr(R), \\ tr(RD) = 2R + n\textrm{ diag}(R).
$$
Moreover 
$$
SCS' = T^{-1}, \\ SRS' = T^{-1}(Ric)T^{-1}
$$
and the lemma claims that 
$$
\Sigma \approx \frac{T^{-1} - \frac{1}{3}T^{-1}(Ric)T^{-1} - \frac{n}{6}S\textrm{diag}(S'(Ric)S)S' } {1 - \frac{1}{6}tr(T^{-1} Ric)}, 
$$
which is different from the approximation $\Sigma \approx T^{-1} - \frac{1}{3}T^{-1}(Ric)T^{-1}$ given in \cite{pennec-riemann}.
\end{densapprox_example_4}
%

{\section{Standard normal distribution on the unit sphere}}

The folded normal distribution on the sphere $\mathbb{S}^n$ is given by 
\begin{equation}\label{eq:normal_density}
dQ(p) = k(2\pi)^{-n/2}\exp(-\frac{1}{2}T(\widetilde{Log_qp}, \widetilde{Log_qp})) dV(p),
\end{equation}
with the following extended expression 
$$
dQ(p) = k(2\pi)^{-n/2}\sum_{i=0}^{\infty} \exp(-\frac{1}{2}(1 \pm \frac{2\pi i}{||Log_qp||})^2T(Log_qp, Log_qp)) dV(p).
$$
Above we sum two terms for each $i$; this is what $\pm$ stands for.

In particular, if we assume that in normal coordinates $v$, $T={\sigma^{-2}}I_n$, 
then the Euclidean standard normal density $k(2\pi)^{-n/2}\exp(-\frac{1}{2}v'Tv)dv$ 
has covariance $C=I_n$ and kurtosis matrix $D=\{D_{kl}^{ij}\}$ 
such that $D_{kl}^{kl} = D_{kl}^{lk} = 1$, for $k\ne l$, 
$D_{kk}^{ii} = 1$, for $k\ne i$, 
$D_{kk}^{kk} = 3$, $k\ne l$ and zero otherwise.

For this particular $T$, the density ($\ref{eq:normal_density}$) is
\begin{equation}\label{eq:stdnormal_density}
dQ(p) = k(2\pi)^{-n/2}\sum_{i=0}^{\infty} \exp(-\frac{1}{2\sigma^2}(||Log_qp|| \pm 2\pi i)^2) dV(p).
\end{equation}
On the sphere, the Ricci tensor matrix is $Ric=I_n$ and 
since $tr(RC) = n\sigma^2$ and $tr(RD) = (n+2)\sigma^2 I_n$
we can simplify (\ref{eq:approx_density_k}) and (\ref{eq:approx_density_cov}) to
$$
k^{-1} = \sigma^n[1-\frac{n}{6}\sigma^2 + O(\sigma^{3})]
$$ and $$
k^{-1}\Sigma = \sigma^n[1 - \frac{n+2}{6}\sigma^2 + O(\sigma^{3})] \sigma^2 I_n.
$$
For n=2, we write
\begin{equation}\label{eq:predict2_sigma_sphere}
\Sigma \approx \frac{1-\frac{2}{3}\sigma^2}{1-\frac{1}{3}\sigma^2} \sigma^2 I_2.
\end{equation}

We can benefit from a better approximation of the volume form and derive more precise estimation than (\ref{eq:predict2_sigma_sphere}).
The volume form of $\mathbb{S}^n$ in normal coordinates $v$ (see for example 2.3 in \cite{chavel-riemannian}) is 
$$
dV(v) = \frac{\sin(||v||)}{||v||}dv
$$
with Taylor expansion
\begin{equation}\label{eq:approx_volume_2sphere}
dV(v) = [1 - \frac{1}{6}||v||^2 + \frac{1}{120}||v||^4 + O(||v||^6)]dv.
\end{equation}
Utilizing the equations 
\begin{enumerate}
\item[(i)] 
$$
\int_{\mathbb{R}^n} (vv')exp(-\frac{1}{2\sigma^2}v'v)dv = \sigma^3 (2\pi)^{n/2} I_n
$$
\item[(ii)] 
$$
\int_{\mathbb{R}^n} (vv')(v'v)exp(-\frac{1}{2\sigma^2}v'v)dv = (n+2)\sigma^5 (2\pi)^{n/2} I_n
$$
\item[(iii)]
$$
\int_{\mathbb{R}^n} (v'v)exp(-\frac{1}{2\sigma^2}v'v)dv = n\sigma^3 (2\pi)^{n/2}.
$$
\item[(iv)] 
$$
\int_{\mathbb{R}^n} (vv')(v'v)^2exp(-\frac{1}{2\sigma^2}v'v)dv = (n^2+3n+11)\sigma^7 (2\pi)^{n/2} I_n
$$
\item[(v)]
$$
\int_{\mathbb{R}^n} (v'v)^2exp(-\frac{1}{2\sigma^2}v'v)dv = n(n+2)\sigma^5 (2\pi)^{n/2}.
$$
\end{enumerate}
one can show following 
\newtheorem{lemma_approx_stdnormal_density}[lemma_approx_density]{Lemma}
\begin{lemma_approx_stdnormal_density}
The standard normal density on $\mathbb{S}^n$ given by (\ref{eq:stdnormal_density}) has
$$
k^{-1} = (2\pi)^{n/2}\sigma^n[1 - \frac{n}{6}\sigma^2 + \frac{n(n+2)}{120}\sigma^4 + O(\sigma^6)],
$$
and 
$$
k^{-1}\Sigma = (2\pi)^{n/2}\sigma^{n+2}[1-\frac{(n+2)}{6}\sigma^2 + \frac{(n^2+3n+11)}{120}\sigma^4 + O(\sigma^6)] I_n.
$$
\end{lemma_approx_stdnormal_density}
%

In particular, for n = 2,
\begin{equation}\label{eq:predict_sigma_sphere}
\Sigma \approx \frac{1-\frac{2}{3}\sigma^2+\frac{7}{40}\sigma^4}{1-\frac{1}{3}\sigma^2+\frac{1}{15}\sigma^4} \sigma^2 I_2,
\end{equation}
and we expect $tr(\hat\Sigma)$ to be underestimate for $\sigma^2$. 

This conclusion we confirm by simulation studies. 
\begin{figure}
\centering
\begin{tabular}{cc}
{\includegraphics[scale=0.4]{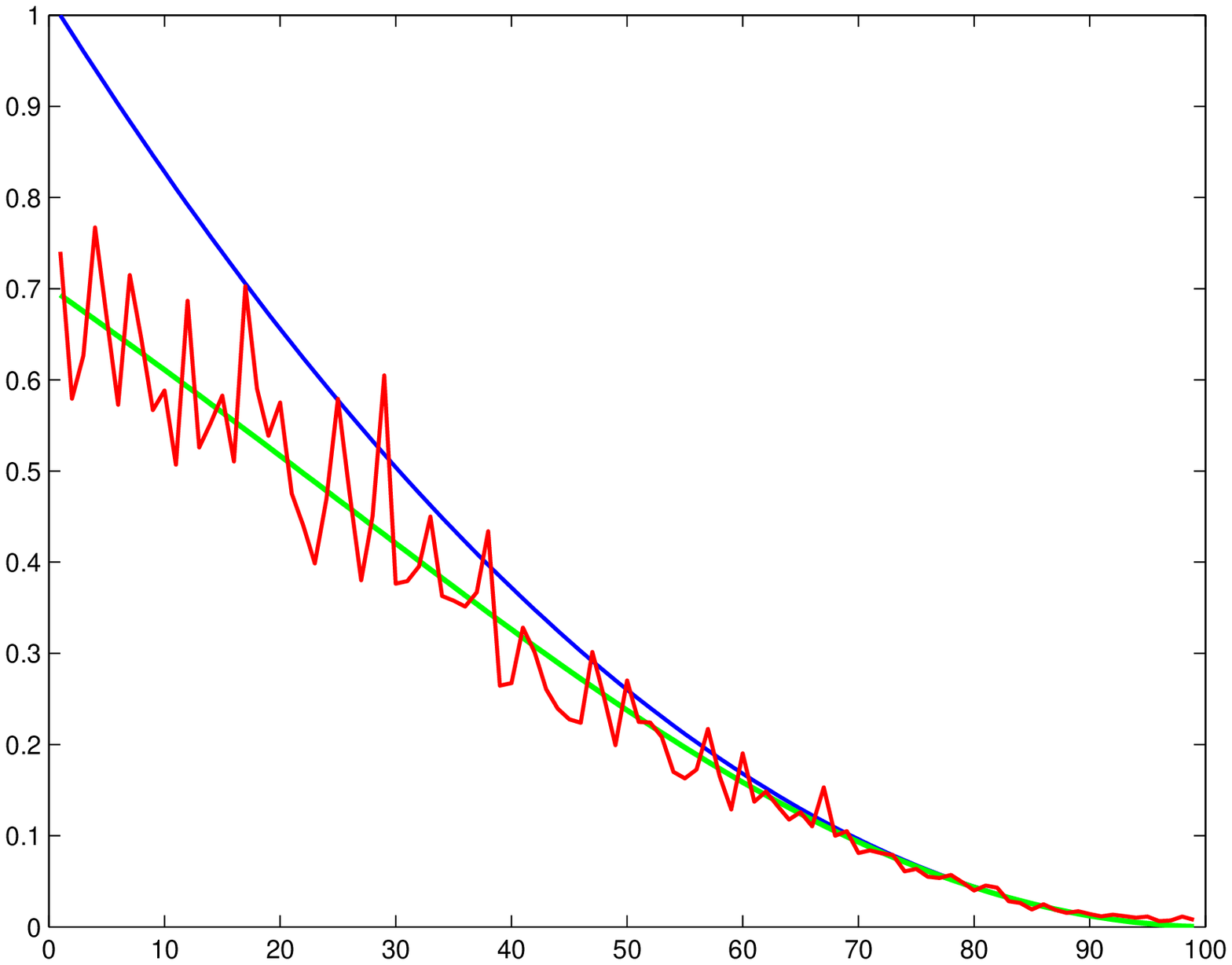}}
{\includegraphics[scale=0.4]{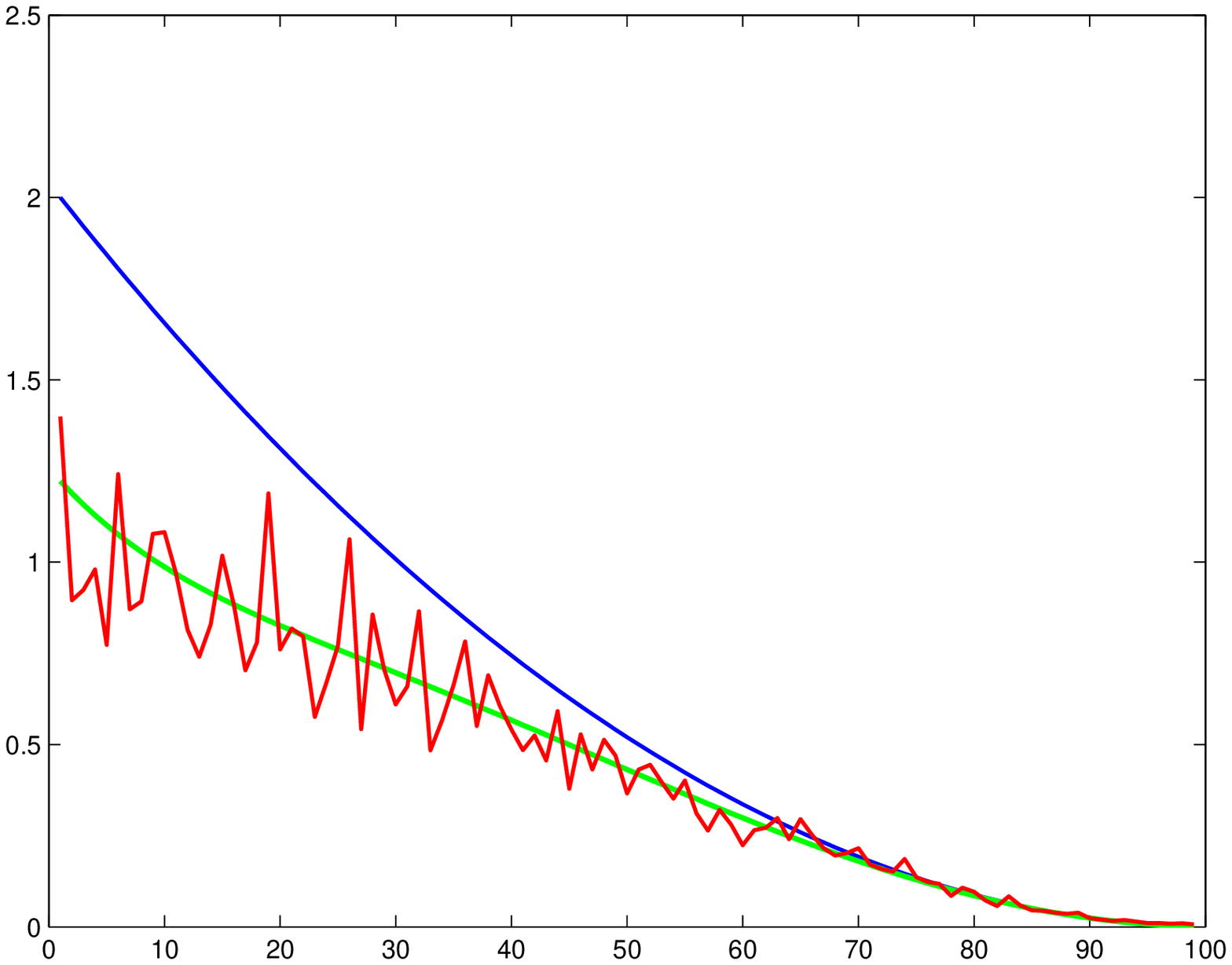}}
\end{tabular}
\caption{Estimation of $\sigma^2$ of normal distribution on $\mathbb{S}^2$ with $T=\sigma^{-2}I_2$. 
Values of $\sigma$ are in blue and decreases from 1 to 0.01 in the left figure and from $\sqrt{2}$ to $\sqrt{2}/100$ in the right one. 
Green curves correspond to the prediction function 
$\frac{1-\frac{2}{3}\sigma^2+\frac{7}{40}\sigma^4}{1-\frac{1}{3}\sigma^2+\frac{1}{15}\sigma^4} \sigma^2$,
as given by equation (\ref{eq:predict_sigma_sphere}). 
Red curves show the estimates $\hat\sigma^2$ calculated using 150 samples for each $\sigma$. 
}
\label{fig:densapprox_sigma_sphere}
\end{figure}
%
Figure (\ref{fig:densapprox_sigma_sphere}) shows the results from our experiment. 
Let (x,y,z) be the cartesian coordinates in $\mathbb{R}^3$.
We generate samples from a normal distribution 
with mean $q=(0,1,0)$ and $T=\sigma^2I_2$ for different values of $\sigma$ shown in blue. 
For every value of $\sigma$, 100 samples are drawn to estimate the covariance $\hat\Sigma$.
The green curve shows the prediction according to (\ref{eq:predict_sigma_sphere}).
The red one shows $\hat\sigma^2 = tr\hat\Sigma$. 
As we see for $ n=2 $ and $ tr(T^{-1})<1 $, $ \hat\sigma^2 $ stays close to the predicted value (\ref{eq:predict_sigma_sphere}).

{\section{Normal distribution on hyperbolic spaces}}

The hyperbolic space $\mathbb{H}^n$ is a Riemannian n-manifold, defined as the half-space 
$\{(x_1,...,x_n), x_n>0\}$ of $\mathbb{R}^n$ endowed with the metric represented by
$$
g_{ij}(x) = \frac{\delta_{ij}}{x_n^2}.
$$
$\mathbb{H}^n$ is geodesically complete and for any point $q\in \mathbb{H}^n0$, 
the exponential map at $q$, $Exp_{q}:\mathbb{R}^n\to \mathbb{H}^n$ is a diffeomorphism on the whole tangent space. 
Thus, the cut locus, $Cut(q)$, is empty. 
It is said that $\mathbb{H}^n$ is a {\it manifold with a pole}.

A normal distribution on $\mathbb{H}^2$ is given by 
\begin{equation}\label{eq:normal_density2}
dQ(p) = k(2\pi)^{-1}\exp(-\frac{1}{2}T(Log_qp, Log_qp)) dV(p).
\end{equation}
In particular, if we assume that in normal coordinates v, $T={\sigma^{-2}}I_n$, 
then 
$$
dQ(v) = k(2\pi)^{-n/2} \exp(-\frac{1}{2\sigma^2}||v||^2) dV(v).
$$
The hyperbolic plane has a constant curvature of -1 and the Ricci tensor matrix is $Ric=-I_n$ 
(for details see \cite{carmo-riemannian}, ch. 8.3).

In two dimensional case, $n=2$, 
we can simplify (\ref{eq:approx_density_k}) and (\ref{eq:approx_density_cov}) to
$$
k^{-1} = \sigma^2[1+\frac{1}{3}\sigma^2 + O(\sigma^4)]
$$$$
k^{-1}\Sigma = \sigma^2[1 + \frac{2}{3}\sigma^2 + O(\sigma^4)] \sigma^2 I_2,
$$
and
\begin{equation}\label{eq:predict2_sigma_hyperplane}
\Sigma \approx \frac{1+\frac{2}{3}\sigma^2}{1+\frac{1}{3}\sigma^2} \sigma^2 I_2.
\end{equation}
We will derive a much better covariance approximation using more precise volume expression.

The volume form of hyperbolic n-manifold $\mathbb{H}^n$ is (see for eample 2.3 in \cite{chavel-riemannian}) 
$$
dV(v) = \frac{sinh(||v||)}{||v||}dv = \frac{\exp(||v||) - \exp(-||v||)}{2||v||}dv, 
$$
and consequently 
\begin{equation}\label{eq:approx_volume_hyperbol}
dV(v) = [1 + \frac{1}{6}||v||^2 + \frac{1}{120}||v||^4 + O(||v||^6)]dv.
\end{equation}

Similarly to the unit n-sphere case, we obtain
\newtheorem{lemma_approx_stdnormal_hyperbol}[lemma_approx_density]{Lemma}
\begin{lemma_approx_stdnormal_hyperbol}
The standard normal density on $\mathbb{H}^n$ has
$$
k^{-1} = (2\pi)^{n/2}\sigma^n[1 + \frac{n}{6}\sigma^2 + \frac{n(n+2)}{120}\sigma^4 + O(\sigma^6)],
$$
and 
$$
k^{-1}\Sigma = (2\pi)^{n/2}\sigma^{n+2}[1+\frac{(n+2)}{6}\sigma^2 + \frac{(n^2+3n+11)}{120}\sigma^4 + O(\sigma^6)] I_n.
$$
\end{lemma_approx_stdnormal_hyperbol}
In particular, for $n=2$
\begin{equation}\label{eq:predict_sigma_hyperbol}
\Sigma \approx \frac{1+\frac{2}{3}\sigma^2+\frac{7}{40}\sigma^4}{1+\frac{1}{3}\sigma^2+\frac{1}{15}\sigma^4} \sigma^2 I_2,
\end{equation}
\begin{figure}
\centering
\begin{tabular}{cc}
{\includegraphics[scale=0.4]{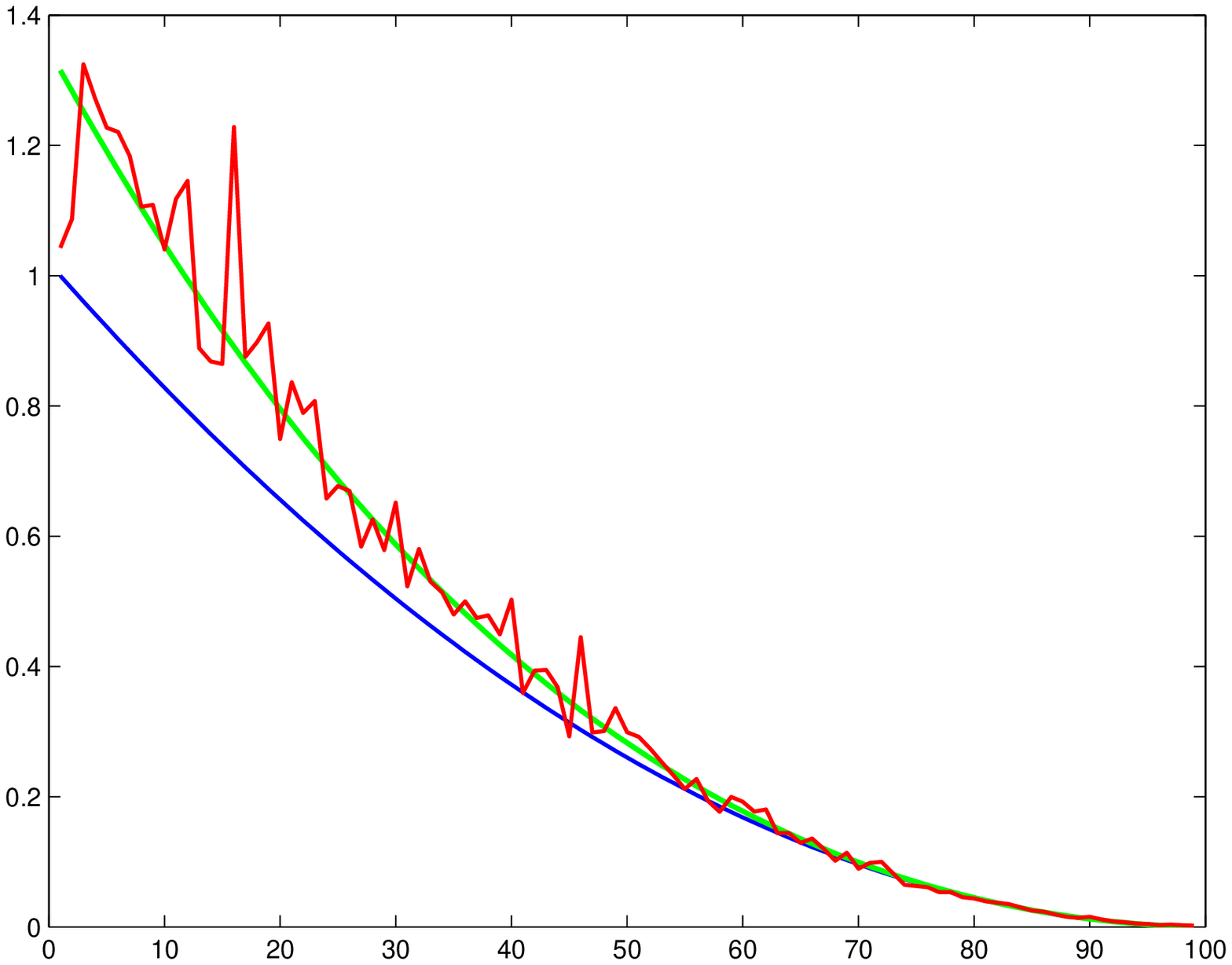}} & 
{\includegraphics[scale=0.4]{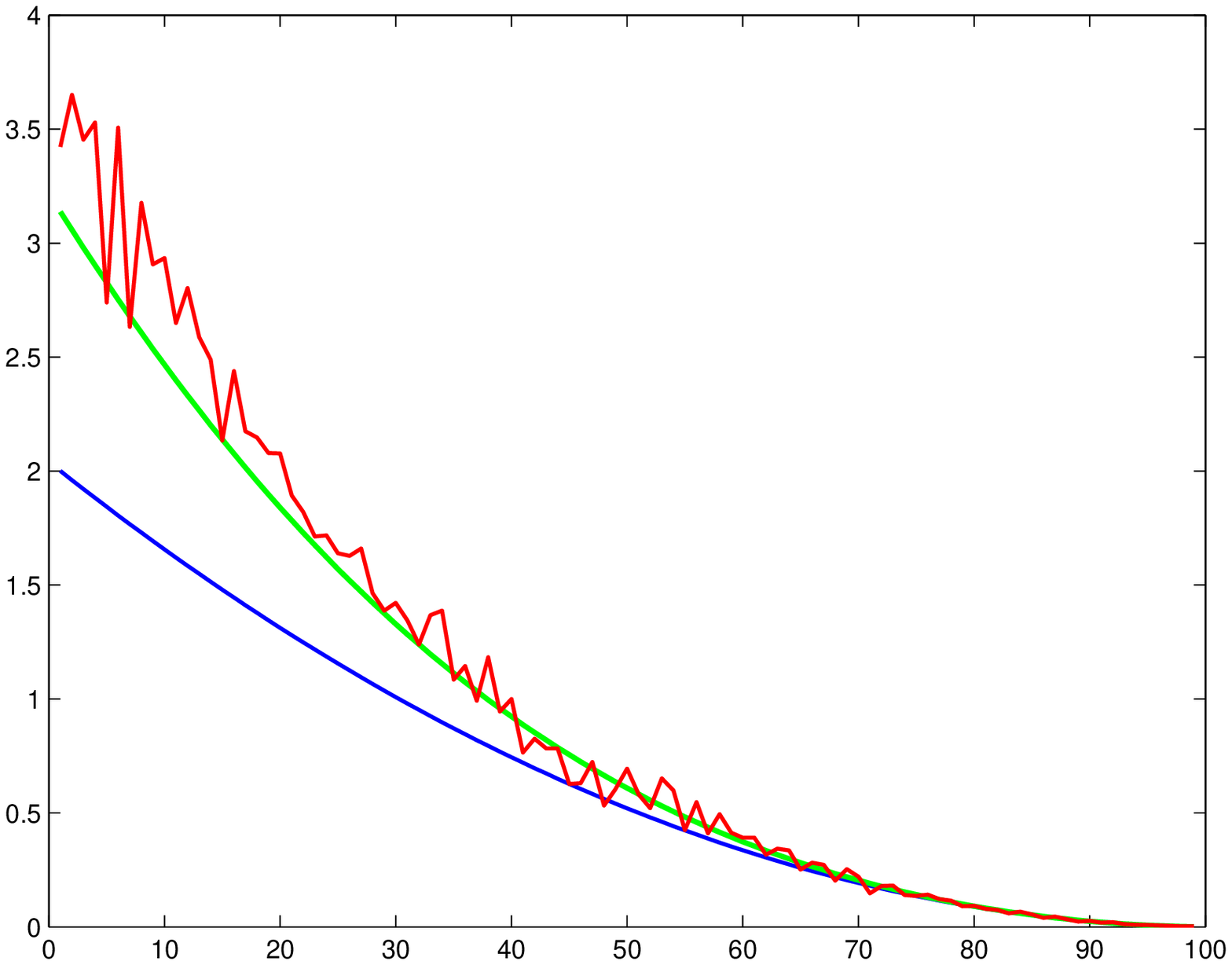}}
\end{tabular}
\caption{Estimation of $\sigma^2$ for normal distribution on $\mathbb{H}^2$ with $T=\sigma^{-2}I_2$. 
Values of $\sigma$ are in blue and decreases from 1 to 0.01 in the left figure and from $\sqrt{2}$ to $\sqrt{2}/100$ in the right one. 
Green curves correspond to the prediction function 
$\frac{1+\frac{2}{3}\sigma^2+\frac{7}{40}\sigma^4}{1+\frac{1}{3}\sigma^2+\frac{1}{15}\sigma^4} \sigma^2$,
as given by equation (\ref{eq:predict_sigma_hyperbol}). 
Red curves show the estimates $\hat\sigma^2$ calculated using 200 samples  for each $\sigma$. 
}
\label{fig:densapprox_sigma_hyperbol}
\end{figure}
Therefore we expect $tr(\hat\Sigma)$ to overestimate $\sigma^2$.
This conlcusion we confirm experimentally (see Figure (\ref{fig:densapprox_sigma_hyperbol})). 
When $\sigma^2<1 $, $ \hat\sigma^2$ stays close to the predicted value (\ref{eq:predict_sigma_hyperbol}). 
For larger values of $\sigma^2$ more precise approximation is needed.

Upon request we provide MATLAB programs for the experiments shown in 
Figures (\ref{fig:densapprox_sigma_sphere}) and (\ref{fig:densapprox_sigma_hyperbol}).

{\section{Summary}}


In this study we try to be more precise and general when defining distributions on 
complete Riemannian manifolds and on compact manifolds in particular. 
We give a consistent definition that accounts for the lack of global parametrization on manifolds  
by being coordinate independent. 
Also, coordinate specific attributes, like concentration matrix and covariance, 
are treated more carefully. They are considered as tensors of appropriate variety. 
The motivating idea behind this point of view is that only coordinate invariant objects should be used for 
statistical inference purposes.

The families of centered distributions we dealt with, 
are usually based on Euclidean multivariate kernel, like the normal one.
That makes the problem of relating the covariance of manifold variable to its Euclidean counterpart interesting.
We expressed formally one possible relation in this regard 
and confirmed it with simulations. Our experiments include normal distribution on the unit 2-sphere, 
which is of interest of directional statistics, and normal distribution on the hyperbolic plane, 
which lack application potential for the moment, but it is an interesting demonstration by itself for 
clearly showing the impact of the negative curvature of the domain.



\end{document}